\documentclass[]{amsart}

\usepackage[latin1]{inputenc}
 \usepackage[english]{babel}
 \usepackage{amsmath}
 \usepackage{amsthm}
 \usepackage{paralist}
 \usepackage{amssymb}
 \usepackage{bibgerm}
 \usepackage{dsfont}
 \usepackage[pass]{geometry}

\usepackage[numbers,sort&compress]{natbib}
\bibpunct[, ]{[}{]}{,}{n}{,}{,}

\theoremstyle{plain}
\newtheorem{theorem}{Theorem}[section]
\newtheorem{lemma}[theorem]{Lemma}

\newtheorem{proposition}[theorem]{Proposition}

\theoremstyle{definition}
\newtheorem{definition}[theorem]{Definition}

\newtheorem{assumption}[theorem]{Assumption}

\theoremstyle{remark}
\newtheorem{remark}{Remark}

\begin{document}
 
\title[Pathwise Uniqueness for SDEs with Singular Drift and Nonconst.~Diffusion~~]{Pathwise Uniqueness for SDEs with Singular Drift and Nonconstant Diffusion:\\A simple proof}
 
\author{Katharina von der L\"uhe}
\address{Katharina von der L\"uhe, Fakult\"at f\"ur Mathematik, Universit\"at Bielefeld, D-33501 Bielefeld, Germany.}
\email{katharina.vonderluehe@math.uni-bielefeld.de}

\begin{abstract}
 A new proof of pathwise uniqueness for SDEs with Sobolev diffusion and integrable drift term is introduced by extending a method from E.\,Fedrizzi and F.\,Flandoli (\cite{FF11}) to the case of nonconstant diffusion. 
\end{abstract}

\maketitle



\section{Introduction}

Let us consider the following stochastic differential equation (SDE):
\[X_t=x+\int\limits_0^tb(s,X_s)ds+\int\limits_0^t\sigma(s,X_s)dW_s,\quad t\in[0,T],\]
where $x\in\mathbb{R}^d$, $b$, $\sigma$ are measurable functions from $[0,T]\times\mathbb{R}^d$ to $\mathbb{R}^d$, respectively $\mathbb{R}^{d\times m}$, and $W$ is an $m$-dimensional standard Wiener process.\\
There are many papers which investigate the problem of existence or uniqueness of solutions for this kind of equation. In addition to the well known result for Lipschitz coefficients by K.\,It\^o, \cite{Ito}, let us mention some of these results here. Strong existence and uniqueness have been obtained for example under local weak monotonicity and weak coercivity conditions on the coefficients. A proof can be found in Chapter 3 of the monograph by W.\,Liu and M.\,R\"ockner \cite{Roeckner_Liu}. Furthermore, in their work \cite{Fang_Zhang}, S.\,Fang and T.\,Zhang relaxed the Lipschitz conditions by a logarithmic factor. Moreover, A.\,Yu.\,Veretennikov proved strong existence and uniqueness for bounded measurable coefficients if the diffusion matrix is nondegenerated, continuous and Lip\-schitz continuous in the spacial variable, see \cite{Veretennikov}. In \cite{Gyongy_Martinez} I.\,Gy\"ongy and T.\,Mart\'inez relaxed this to locally unbounded drifts, namely $b\in L_{loc}^{2(d+1)}(\mathbb{R}_+\times\mathbb{R}^d)$ and $b$ almost everywhere bounded by a constant plus some nonnegative function in $L^{d+1}(\mathbb{R}_+\times\mathbb{R}^d)$.\\
In \cite{KR05} N.\,Krylov and M.\,R\"ockner proved the existence of a unique strong solution up to some explosion time in the case where the diffusion coefficient $\sigma$ is the unit matrix and the drift coefficient $b$ is in $L_{\text{loc}}^q(\mathbb{R}_+;L_{\text{loc}}^p(\mathbb{R}^d))$ for some $p,q>1$ fulfilling
\begin{equation}\label{Prodi_Serrin}
 \frac{d}{p}+\frac{2}{q}<1.
\end{equation}\\
If the diffusion is not constant and nondegenerate it is also possible to get strong existence and uniqueness results under similar conditions on the drift. The most general result can be found in the work of X.\,Zhang \cite{Zhang_Neu}, respectively for the case $p=q$ see \cite{Zhang}. There, the drift is again in $L^q_\text{loc}(\mathbb{R}_+,L^p(\mathbb{R}^d))$ for $p,q>1$ fulfilling \eqref{Prodi_Serrin}. The diffusion coefficient is uniformly continuous in space, locally uniformly with respect to time, nondegenerated, bounded and the gradient is also in 
$L_{loc}^q(\mathbb{R}_+,L^p(\mathbb{R}^d))$. The idea of the proof is to remove the drift by the so-called Zvonkin transformation, see \cite{Zvonkin}, and use known results for SDEs with zero drift. This transformation is based on the solution $u$ to the equation
\[\partial_tu+\sum\limits_{i=1}^db^i\partial_{x_i}u+\frac{1}{2}\sum\limits_{i=1}^d\sum\limits_{j=1}^d(\sigma\sigma^*)_{ij}\partial_{x_ix_j}^2u=0, ~~~~~~u(T,x)=x.\]
Then one gets a one-to-one correspondence between solutions $X_t$ for the original SDE and solutions $u(t,X_t)$ for the transformed equation without drift term.\\
In the case of constant $\sigma$ there is a much simpler proof for the pathwise uniqueness which is due to E.\,Fedrizzi and F.\,Flandoli (see \cite{FF11}) under similar conditions as in \cite{KR05}. They gave an elementary and short proof by developing another transformation of the SDE. The aim of this work is to extend their method to include the case of $b$ and $\sigma$ under the conditions as in \cite{Zhang_Neu}. The nonconstant diffusion leads to additional terms when performing the transformation of \cite{FF11} which have to be controlled. One of the main tools to overcome these difficulties are Krylov Estimates but the price to pay is that we have to assume that \eqref{Prodi_Serrin} holds with $1/2$ replacing $1$ on its right hand side.\\
For simplicity we will state our result under global assumptions, but there are no difficulties to extend it by localization techniques, e.g.\,in the same way as in \cite{Zhang_Neu}.\\

\section{Preliminaries and main result}

\begin{definition}
 For $p,q\in(1,\infty)$ we define
 \[\| f\|_{L_p^q(T)}~:=~\left(\int\limits_0^T\left(~\int\limits_{\mathbb{R}^d}|f(t,x)|^p\,dx\right)^\frac{q}{p}\,dt\right)^\frac{1}{q},\]
 where $|\cdot|$ denotes the Hilbert--Schmidt norm.\\
 We define $L_p^q(T)$ to be the space of measurable functions $f:[0,T]\times\mathbb{R}^d\to\mathbb{R}^d$ (respectively $\mathbb{R}^{d\times m}$) such that $\|f\|_{L_p^q(T)}<\infty$.\\
 Furthermore, we set
 \[W_{q,p}^{1,2}(T):=\left\{f:[0,T]\times\mathbb{R}^d\to\mathbb{R}^d~\Big|~f,\partial_tf,\partial_xf,\partial_x^2f\in L_p^q(T)\right\},\]
 where $\partial_t$, $\partial_x$, $\partial_x^2$ denote the weak derivatives with respect to time, respectively space. The associated norm is given by
 \[\| f\|_{W_{q,p}^{1,2}(T)}~:=~\| f\|_{L_p^q(T)}+\|\partial_tf\|_{L_p^q(T)}+\|\partial_xf\|_{L_p^q(T)}+\|\partial_x^2f\|_{L_p^q(T)}.\]
\end{definition}
 
We consider the SDE
 \begin{equation}\label{SDE}
 X_t =x+\int\limits_0^tb(s,X_s)\,ds+\int\limits_0^t\sigma(s,X_s)\,dW_s, \quad t\in[0,T],
 \end{equation}
where $W$ is an $m$-dimensional standard Wiener process on a filtered probability space $(\Omega,(\mathcal{F}_t)_t,\mathbb{P})$, with $(\mathcal{F}_t)_t$ fulfilling the usual conditions, $x\in\mathbb{R}^d$ and $b:[0,T]\times\mathbb{R}^d\to\mathbb{R}^d$, $\sigma:[0,T]\times\mathbb{R}^d\to\mathbb{R}^{d\times m}$ are measurable functions with the following properties:

\begin{assumption}\label{assumptions}
For some $p,q>1$ with
\[\frac{d}{p}+\frac{2}{q}<\frac{1}{2},\] 
we have
\begin{enumerate}
\item[(c1)] $b\in L_p^q(T)$,
\item[(c2)] $\sigma$ is uniformly continuous in $x$, uniformly with respect to $t$, i.e. for all $\varepsilon>0$ exists a $\delta>0$ such that
\[\sup_{t\in[0,T]}|\sigma(t,x)-\sigma(t,y)|<\varepsilon\quad\text{for all}~x,y\in\mathbb{R}^d~\text{with}~|x-y|<\delta.\]
\item[(c3)] $\sigma$ is nondegenerated, i.e.\,there exists a constant $c_\sigma>0$ such that
  \[\langle\sigma\sigma^*(t,x)\xi,\xi\rangle\geq c_\sigma\langle I\xi,\xi\rangle~~~~~~\forall~\xi\in\mathbb{R}^d~~~\forall~(t,x)\in[0,T]\times\mathbb{R}^d,\]
  where $\sigma^*$ denotes the transposed matrix of $\sigma$,
\item[(c4)] $\sigma$ is bounded by a constant $\tilde{c}_\sigma$,
\item[(c5)] $\partial_x\sigma\in L_p^q(T)$.
\end{enumerate}
\end{assumption}

\begin{definition}(weak/strong solution)
 A weak solution for equation $(\ref{SDE})$ is a pair $(X,W)$ on a filtered probability space $(\Omega,(\mathcal{F}_t)_t,\mathbb{P})$ such that $X$ is continuous, $(\mathcal{F}_t)_t$-adapted, 
 fulfills
 \begin{equation}\label{endliche variation}
  \mathbb{P}\left(\int\limits_0^T|b(s,X_s)|\,ds<\infty\right)=1,
 \end{equation}
 \begin{equation}\label{existence stochastic integral}
  \mathbb{P}\left(\int\limits_0^T|\sigma(s,X_s)|^2\,ds<\infty\right)=1,
 \end{equation}
 $W_t$ is an $\mathcal{F}_t$-Brownian motion and $(X,W)$ satisfies equation $(\ref{SDE})$ almost surely.\\
 Given a Brownian motion $W$ on a probability space, a strong solution for equation $(\ref{SDE})$ is a continuous process $X$ which is adapted to the filtration generated by $W$, fulfills \eqref{endliche variation}, \eqref{existence stochastic integral} and satisfies equation $(\ref{SDE})$ almost surely.
\end{definition}

\begin{definition}(Pathwise Uniqueness)
 We say that pathwise uniqueness holds for equation $(\ref{SDE})$ if for two weak solutions $(X,W)$, $(\tilde{X},\tilde{W})$, defined on the same probability space, we have that $X_0=\tilde{X}_0$ and $W=\tilde{W}$ imply
 \[\mathbb{P}\left(X_t=\tilde{X}_t~~~\forall~t\in[0,T]\right)=1.\]
\end{definition}

\begin{theorem}\label{main_result}
 Under Assumption \ref{assumptions}, we have pathwise uniqueness for equation \eqref{SDE}.
 
\end{theorem}

\begin{remark}
One obtains the same result for $p,q>2(d+1)$ if instead of Assumption \ref{assumptions} $(c2)$ $\sigma$ is assumed to be continuous and such that for all $f\in L_p^{q}(T)$ there exists a solution to the equation
\[\partial_tu+\frac{1}{2}\sum\limits_{i=1}^d\sum\limits_{j=1}^d(\sigma\sigma^*)_{ij}\partial_{x_ix_j}^2u=f \quad\text{on}~[0,T],\quad u(T,x)=0,\]
 such that
 \[\|u\|_{L_p^q(T)}\leq C\|f\|_{L_p^q(T)},\]
 where $C$ is independent of $f$ and increasing in $T$. For details see \cite{vdl}.
\end{remark}

In the following, whenever we speak of two solutions, we mean two weak solutions defined on the same probability space with the same Brownian motion.\\
 Furthermore by $C>0$ we always denote various finite constants, where we often indicate the dependence of parameters by writing them in brackets.

\section{Transformation of the SDE}
The following transformation works analogously to the transformation of \cite{FF11} despite the appearance of additional terms in the partial differential equations and the stochastic integrals.\\
Assume that $b$ and $\sigma$ fulfill Assumption \ref{assumptions}. Then by Theorem 10.3 and Remarks 10.4 and 10.5 in \cite{KR05} for every $f\in L_p^q(T)$ there exists a solution $u\in W_{q,p}^{1,2}(T)$ to the equation
\begin{equation}\label{PDE_f}
 \partial_tu+\frac{1}{2}\sum\limits_{i=1}^d\sum\limits_{j=1}^d(\sigma\sigma^*)_{ij}\partial_{x_ix_j}^2u=f\quad\text{on}~[0,T], \quad u(T,x)=0
\end{equation}
such that
\begin{equation}\label{abschaetzung_loesung_pde}
 \| u\|_{W_{q,p}^{1,2}(T)}\leq C\| f\|_{L_p^q(T)},
\end{equation}
where $C$ does not depend on $f$ and is increasing in $T$.
Then by the H\"older continuity of $\partial_xu$, see \cite{KR05} Lemma 10.2, we have
\begin{equation}\label{abschaetzung_fuer_erste_ableitung}
 \sup_{(t,x)\in[0,T]\times\mathbb{R}^d}|\partial_xu(t,x)| \leq C(p,q,\varepsilon,T)T^\frac{\varepsilon}{2}\| f\|_{L_p^q(T)}
\end{equation}
for every $\varepsilon\in(0,1)$, which fulfills
\[\varepsilon+\frac{d}{p}+\frac{2}{q}<1,\]
with $C(p,q,\varepsilon,T)$ increasing in $T$. We can therefore assume the constant in front of $\|f\|_{L_p^q(T)}$ to be as small as we want by choosing $T$ appropriate which will be of importance in Lemma \ref{Ungleichungen}.
Now, let $U_b$ a solution to the equation
\begin{equation}\label{PDE_-b}
\partial_tu+\frac{1}{2}\sum\limits_{i=1}^d\sum\limits_{j=1}^d(\sigma\sigma^*)_{ij}\partial_{x_ix_j}^2u=-b \quad\text{on}~[0,T],\quad u(T,x)=0.
\end{equation}
Using It\^o's formula for functions in $W_{q,p}^{1,2}(T)$ (Proposition \ref{Ito-Formel}) and that $U_b$ is a solution to PDE \eqref{PDE_-b}, we get
\begin{align*}
U_b(t,X_t) & = U_b(0,x) + \int\limits_0^t\partial_xU_b(s,X_s)b(s,X_s)\,ds + \int\limits_0^t\partial_xU_b(s,X_s)\sigma(s,X_s)\,dW_s\\
           & ~~~ - \int\limits_0^tb(s,X_s)\,ds.
\end{align*}
That implies
\begin{align*}
\int\limits_0^tb(s,X_s)\,ds & = U_b(0,x)-U_b(t,X_t) + \int\limits_0^t\partial_xU_b(s,X_s)b(s,X_s)\,ds\\
                    & ~~~~~~~ + \int\limits_0^t\partial_xU_b(s,X_s)\sigma(s,X_s)\,dW_s.
\end{align*}
Now, we define
\[\mathcal{T}(b):=\partial_xU_b\cdot b\]
and transform SDE \eqref{SDE} by replacing the drift term:
\begin{align}
X_t & =x+U_b(0,x)-U_b(t,X_t)+\int\limits_0^t\mathcal{T}(b)(s,X_s)\,ds\notag\\
    & ~~~~~~~~~~~+\int\limits_0^t\partial_xU_b(s,X_s)\sigma(s,X_s)+\sigma(s,X_s)\,dW_s.\label{SDE_transformed1}
\end{align}
Note, that $\mathcal{T}(b)\in L_p^q(T)$ since $\partial_xU_b$ is bounded and $b\in L_p^q(T)$. Next, let $U_{\mathcal{T}(b)}$ be a solution to the equation
\[\partial_tu+\frac{1}{2}\sum\limits_{i=1}^d\sum\limits_{j=1}^d(\sigma\sigma^*)_{ij}\partial_{x_ix_j}^2u=-\mathcal{T}(b)\quad\text{on}~[0,T],\quad u(T,x)=0.\]
Using again It\^o's formula (Proposition \ref{Ito-Formel}) and that $U_{\mathcal{T}(b)}$ solves the equation above, we get
\begin{align*}
U_{\mathcal{T}(b)}(t,X_t) & = U_{\mathcal{T}(b)}(0,x) + \int\limits_0^t\partial_xU_{\mathcal{T}(b)}(s,X_s)b(s,X_s)\,ds\\
                          & ~~~ + \int\limits_0^t\partial_xU_{\mathcal{T}(b)}(s,X_s)\sigma(s,X_s)\,dW_s - \int\limits_0^t\mathcal{T}(b)(s,X_s)\,ds,
\end{align*}
and therefore
\begin{align*}
\int\limits_0^t\mathcal{T}(b)(s,X_s)\,ds & = U_{\mathcal{T}(b)}(0,x) - U_{\mathcal{T}(b)}(t,X_t) + \int\limits_0^t\partial_xU_{\mathcal{T}(b)}(s,X_s)b(s,X_s)\,ds\\
                                         & ~~~~~~~ + \int\limits_0^t\partial_xU_{\mathcal{T}(b)}(s,X_s)\sigma(s,X_s)\,dW_s.
\end{align*}
As before, we define
\[\mathcal{T}^2(b):=\partial_xU_{\mathcal{T}(b)}\cdot b\]
and replace the drift term in the transformed SDE \eqref{SDE_transformed1}:
\begin{align*}
X_t & = x + U_b(0,x) + U_{\mathcal{T}(b)}(0,x) - U_b(t,X_t) - U_{\mathcal{T}(b)}(t,X_t) + \int\limits_0^t\mathcal{T}^2(b)(s,X_s)\,ds\\
    & ~~~ + \int\limits_0^t\partial_xU_b(s,X_s)\sigma(s,X_s)+\partial_xU_{\mathcal{T}(b)}(s,X_s)\sigma(s,X_s)+\sigma(s,X_s)\,dW_s.
\end{align*}
Iteration yields after $n+1$ steps
\begin{align}
X_t & = x + \sum\limits_{k=0}^nU_{\mathcal{T}^k(b)}(0,x) - \sum\limits_{k=0}^nU_{\mathcal{T}^k(b)}(t,X_t) + \int\limits_0^t\mathcal{T}^{n+1}(b)(s,X_s)\,ds\notag\\
    & ~~~ + \int\limits_0^t\sum\limits_{k=0}^n\partial_xU_{\mathcal{T}^k(b)}(s,X_s)\sigma(s,X_s)+\sigma(s,X_s)\,dW_s\label{SDE_transformedn}
\end{align}
with the convention
\[\mathcal{T}^0(b)=b ~~~~\text{and}~~~~\mathcal{T}^{k+1}(b)=\partial_xU_{\mathcal{T}^k(b)}\cdot b.\]
We define
\[U^{(n)}(t,x):=\sum\limits_{k=0}^nU_{\mathcal{T}^k(b)}(t,x)\]
and therefore, SDE \eqref{SDE_transformedn} becomes
\begin{align}
X_t & = x + U^{(n)}(0,x) - U^{(n)}(t,X_t) + \int\limits_0^t\mathcal{T}^{n+1}(b)(s,X_s)\,ds\notag\\
    & ~~~ + \int\limits_0^t\left(\partial_xU^{(n)}(s,X_s)+I\right)\sigma(s,X_s)\,dW_s.\label{SDE_transformedn_neu}
\end{align}
For two solutions $X_t^{(1)}$, $X_t^{(2)}$ we define
\begin{align*}
              Y_t^{(i,n)} & := X_t^{(i)}+U^{(n)}(t,X_t^{(i)}),\quad i=1,2,~\text{and}\\
     b^{(n)}(t,x) & := \mathcal{T}^{n+1}(b)(t,x),\\
\sigma^{(n)}(t,x) & := \left(\partial_xU^{(n)}(t,x)+I\right)\sigma(t,x).
\end{align*}
Then equation \eqref{SDE_transformedn_neu} reads
\begin{equation}\label{neue_formulierung_SDE}
Y_t^{(i,n)} = Y_0^{(i,n)} + \int\limits_0^tb^{(n)}(s,X_s^{(i)})\,ds + \int\limits_0^t\sigma^{(n)}(s,X_s^{(i)})\,dW_s.
\end{equation}

The following Lemma summarizes some properties of the transformed equation which are necessary in the proof of pathwise uniqueness. It is similar to Lemma 7 in \cite{FF11} and so is the proof.
\begin{lemma}\label{Ungleichungen}
 Let $(c1)-(c4)$ of Assumption \ref{assumptions} be fulfilled and $X_t^{(1)}$, $X_t^{(2)}$ be two solutions to \eqref{SDE}. Then there exists $0<T_0\leq T$ such that for all $T'\in(0,T_0]$ we have 
 \begin{enumerate}
 \item[\normalfont(i)] $\begin{aligned}[t]\|b^{(n)}\|_{L_p^q(T')}\leq\frac{1}{2^{n+1}}\|b\|_{L_p^q(T')},\end{aligned}$
 \item[\normalfont(ii)] $\begin{aligned}[t]\sum\limits_{k=0}^n\sup_{(t,x)\in[0,T']\times\mathbb{R}^d}|\partial_xU_{\mathcal{T}^k(b)}(t,x)|\leq\frac{1}{2},\end{aligned}$
 \item[\normalfont(iii)] $\begin{aligned}[t]\|\partial_x^2U^{(n)}\|_{L_p^q(T')}\leq C\end{aligned}$
       for some constant $C>0$, independent of $n$, and
 \item[\normalfont(iv)] $\begin{aligned}[t]\left|Y_t^{(1,n)}-Y_t^{(2,n)}\right| & \leq \frac{3}{2}\left|X_t^{(1)}-X_t^{(2)}\right|,\\
          \left|X_t^{(1)}-X_t^{(2)}\right| & \leq 2\left|Y_t^{(1,n)}-Y_t^{(2,n)}\right|~~\text{for all }t\in(0,T'].\end{aligned}$
 \end{enumerate}
\end{lemma}

\section{Pathwise uniqueness}

We now prove Theorem \ref{main_result}. It works analogously to \cite{FF11}, based on Lemma \ref{Ungleichungen} and three results, namely Lemmas \ref{moments}, \ref{convergence of the drift} and Proposition \ref{uniform bound for A_T}, which are similar to \cite{FF11} but with different proofs. This is due to the fact that in our framework the solution is in general not a Brownian motion. For reasons of readability we defer the proofs to the next section. 
\begin{proof}(Proof of Theorem \ref{main_result} for small $T$)
 In the following, we denote by $x^i$ the $i$-th entry of a vector $x\in\mathbb{R}^d$. Let Assumption \ref{assumptions} be fulfilled and $X_t^{(1)}$, $X_t^{(2)}$ be two solutions to \eqref{SDE}. Furthermore, let $T:=T_0$ from Lemma \ref{Ungleichungen} and $Y_t^{(i,n)}$ given by \eqref{neue_formulierung_SDE}.
 By It\^o's formula and an application of the inequality of Cauchy and Schwarz we then have
 \begin{align}
  d\left|Y_t^{(1,n)}-Y_t^{(2,n)}\right|^2 & \leq 2\left|Y_t^{(1,n)}-Y_t^{(2,n)}\right|\left|b^{(n)}(t,X_t^{(1)})-b^{(n)}(t,X_t^{(2)})\right|\,dt\notag\\
                                          & ~~~ + 2\left\langle Y_t^{(1,n)}-Y_t^{(2,n)},\left(\sigma^{(n)}(t,X_t^{(1)})-\sigma^{(n)}(t,X_t^{(2)})\right)\,dW_t\right\rangle\notag\\
                                          & ~~~ + \left|\sigma^{(n)}(t,X_t^{(1)})-\sigma^{(n)}(t,X_t^{(2)})\right|^2\,dt.\label{Abschaetzung_d_Differenz}
 \end{align}
 Moreover, with
 \[A_t^{(n)}:=\int\limits_0^t\frac{\left|\sigma^{(n)}(s,X_s^{(1)})-\sigma^{(n)}(s,X_s^{(2)})\right|^2}{\left|Y_s^{(1,n)}-Y_s^{(2,n)}\right|^2}\textbf{1}_{\{Y_s^{(1,n)}\neq Y_s^{(2,n)}\}}ds,\]
 we have
 \begin{align*}
  d\left(e^{-A_t^{(n)}}\left|Y_t^{(1,n)}-Y_t^{(2,n)}\right|^2\right) & =  e^{-A_t^{(n)}}d\left|Y_t^{(1,n)}-Y_t^{(2,n)}\right|^2\\
                                                                     & ~~~ - \left|Y_t^{(1,n)}-Y_t^{(2,n)}\right|^2e^{-A_t^{(n)}}dA_t^{(n)},
 \end{align*}
 since the quadratic covariation is zero due to the monotonicity of $e^{-A_t^{(n)}}$.
 Now, we use inequality \eqref{Abschaetzung_d_Differenz} to conclude that
 \begin{align*}
  d & \left(e^{-A_t^{(n)}}\left|Y_t^{(1,n)}-Y_t^{(2,n)}\right|^2\right)\\
    & \leq 2e^{-A_t^{(n)}}\left|Y_t^{(1,n)}-Y_t^{(2,n)}\right|\left|b^{(n)}(t,X_t^{(1)})-b^{(n)}(t,X_t^{(2)})\right|\,dt\\
    & ~~~ + 2e^{-A_t^{(n)}}\left\langle Y_t^{(1,n)}-Y_t^{(2,n)},\left(\sigma^{(n)}(t,X_t^{(1)})-\sigma^{(n)}(t,X_t^{(2)})\right)\,dW_t\right\rangle
 \end{align*}
 and thus,
 \begin{align*}
  \mathbb{E} & \left[e^{-A_t^{(n)}}\left|Y_t^{(1,n)}-Y_t^{(2,n)}\right|^2\right]\\
             & \leq \mathbb{E}\left[\left|Y_0^{(1,n)}-Y_0^{(2,n)}\right|^2\right]\\
             & ~~~ + 2\mathbb{E}\left[\int\limits_0^te^{-A_s^{(n)}}\left|Y_s^{(1,n)}-Y_s^{(2,n)}\right|\left|b^{(n)}(s,X_s^{(1)})-b^{(n)}(s,X_s^{(2)})\right|\,ds\right]\\
             & ~~~ + 2\mathbb{E}\left[\int\limits_0^te^{-A_s^{(n)}}\left\langle Y_s^{(1,n)}-Y_s^{(2,n)},\left(\sigma^{(n)}(s,X_s^{(1)})-\sigma^{(n)}(s,X_s^{(2)})\right)\,dW_s\right\rangle\right].
 \end{align*}
 With the help of Lemma \ref{Ungleichungen}, we get
 \begin{align}
  \mathbb{E} & \left[e^{-A_t^{(n)}}\left|Y_t^{(1,n)}-Y_t^{(2,n)}\right|^2\right]\notag\\
             & \leq \frac{9}{4}\left|x^{(1)}-x^{(2)}\right|^2\notag\\
             & ~~~ + 3\mathbb{E}\left[\int\limits_0^t\left|X_s^{(1)}-X_s^{(2)}\right|\left|b^{(n)}(s,X_s^{(1)})-b^{(n)}(s,X_s^{(2)})\right|\,ds\right]\label{Abschaetzung_Erwartungswert_d_Differenz}\\
             & ~~~ + 2\mathbb{E}\left[\int\limits_0^te^{-A_s^{(n)}}\left\langle Y_s^{(1,n)}-Y_s^{(2,n)},\left(\sigma^{(n)}(s,X_s^{(1)})-\sigma^{(n)}(s,X_s^{(2)})\right)\,dW_s\right\rangle\right].\notag
 \end{align}
 Summarizing, for two solutions with the same initial values, we have for all $t\leq T$
 \begin{align*}
  \mathbb{E}\left[\left|X_{t}^{(1)}-X_{t}^{(2)}\right|\right] & \leq \mathbb{E}\left[2\left|Y_{t}^{(1,n)}-Y_{t}^{(2,n)}\right|\right]\\
                                                              & = 2\mathbb{E}\left[e^{\frac{1}{2}A_{t}^{(n)}}e^{-\frac{1}{2}A_{t}^{(n)}}\left|Y_{t}^{(1,n)}-Y_{t}^{(2,n)}\right|\right]\\
                                                              & \leq 2\mathbb{E}\left[e^{A_{t}^{(n)}}\right]^\frac{1}{2}\mathbb{E}\left[e^{-A_{t}^{(n)}}\left|Y_{t}^{(1,n)}-Y_{t}^{(2,n)}\right|^2\right]^\frac{1}{2}.
 \end{align*}
 With inequality \eqref{Abschaetzung_Erwartungswert_d_Differenz} we obtain
 \begin{align}
  \mathbb{E} & \left[\left|X_{t}^{(1)}-X_{t}^{(2)}\right|\right]\notag\\
             & \leq 2\mathbb{E}\left[e^{A_{t}^{(n)}}\right]^\frac{1}{2}\left(3\mathbb{E}\left[\int\limits_0^T\left|X_s^{(1)}-X_s^{(2)}\right|\left|b^{(n)}(s,X_s^{(1)})-b^{(n)}(s,X_s^{(2)})\right|\,ds\right]\right.\label{Vor_Martingal_Eigenschaft}\\
             & ~~~~~ + \left.2\mathbb{E}\left[\int\limits_0^{t}e^{-A_s^{(n)}}\left\langle Y_s^{(1,n)}-Y_s^{(2,n)},\left(\sigma^{(n)}(s,X_s^{(1)})-\sigma^{(n)}(s,X_s^{(2)})\right)\,dW_s\right\rangle\right]\right)^\frac{1}{2}.\notag
 \end{align}
 Note that the second expectation term vanishes due to the martingale property of the stochastic integral which is well defined as $\sigma^{(n)}$ is bounded and $|Y_t^{(1,n)}-Y_t^{(2,n)}|^2$ is integrable by the following Lemma.
\begin{lemma}\label{moments}
  Let $(c1)-(c4)$ of Assumption \ref{assumptions} be fulfilled.
  If $X_t$ is a solution to SDE \eqref{SDE}, we have
  \[\mathbb{E}\left[\sup_{t\in[0,T]}|X_t|\right]<\infty~~~\text{and}~~~\sup_{t\in[0,T]}\mathbb{E}\left[|X_t|^2\right]<\infty.\]
 \end{lemma}
 Therefore, by \eqref{Vor_Martingal_Eigenschaft} we have
 \begin{align*}
  \mathbb{E} & \left[\left|X_{t}^{(1)}-X_{t}^{(2)}\right|\right]\\
             & \leq C\mathbb{E}\left[e^{A_{t}^{(n)}}\right]^\frac{1}{2}\mathbb{E}\left[\int\limits_0^T\left|X_s^{(1)}-X_s^{(2)}\right|\left|b^{(n)}(s,X_s^{(1)})-b^{(n)}(s,X_s^{(2)})\right|\,ds\right]^\frac{1}{2}\\
             & \leq C\mathbb{E}\left[e^{A_{t}^{(n)}}\right]^\frac{1}{2}\mathbb{E}\left[\int\limits_0^T\left|X_s^{(1)}-X_s^{(2)}\right|^2\,ds\right]^\frac{1}{4}\mathbb{E}\left[\int\limits_0^T\left|b^{(n)}(s,X_s^{(1)})-b^{(n)}(s,X_s^{(2)})\right|^2\,ds\right]^\frac{1}{4}\\
             & \leq C\mathbb{E}\left[e^{A_{t}^{(n)}}\right]^\frac{1}{2}\mathbb{E}\left[\int\limits_0^T\left|b^{(n)}(s,X_s^{(1)})-b^{(n)}(s,X_s^{(2)})\right|^2\,ds\right]^\frac{1}{4}
 \end{align*}
 for all $n\in\mathbb{N}$, where the last inequality follows from Lemma \ref{moments}. The proof of pathwise uniqueness is complete if we show that the first term is uniformly bounded in $n$ and that the second term converges to zero. These assertions are given by the next two statements which we also prove in the next section.
 \begin{lemma}\label{convergence of the drift}
  Let $(c1)-(c4)$ of Assumption \ref{assumptions} be fulfilled and $X_t^{(1)}$, $X_t^{(2)}$ be two solutions of \eqref{SDE}. Then we have
  \[\lim_{n\to\infty}\mathbb{E}\left[\int\limits_0^T\left|b^{(n)}(t,X_t^{(1)})-b^{(n)}(t,X_t^{(2)})\right|^2\,dt\right]=0.\]
 \end{lemma}
 \begin{proposition}\label{uniform bound for A_T}
  Let Assumption \ref{assumptions} be fulfilled and $X_t^{(1)}$, $X_t^{(2)}$ be two solutions to \eqref{SDE}. Then there exists a constant $C>0$ such that
  \[\mathbb{E}\left[e^{A_T^{(n)}}\right]\leq C~~~~~\text{uniformly for all }n\in\mathbb{N}.\]
 \end{proposition}
 Hence, we proved
 \[X_t^{(1)}=X_t^{(2)}~~\mathbb{P}\text{-a.\,s.}~~~\forall~t\in[0,T].\]
 Thus,
 \[\mathbb{P}\left(X_t^{(1)}=X_t^{(2)}~\forall~t\in\mathbb{Q}\cap[0,T]\right)=1\]
 and by continuity of the solutions we obtain
 \[\mathbb{P}\left(X_t^{(1)}=X_t^{(2)}~\forall~t\in[0,T]\right)=1.\]
\end{proof}

\begin{remark}
 The interval of pathwise uniqueness can easily be extended to arbitrarily large $T$ by means of a time-shift argument.
\end{remark}

\section{Proofs of auxiliaries}

\begin{proof}(Proof of Lemma \ref{moments})
 We have
 \[\mathbb{E}\left[\sup_{t\in[0,T]}|X_t|\right]\leq |x| + \mathbb{E}\left[~\int\limits_0^T|b(s,X_s)|\,ds~\right] + \mathbb{E}\left[\sup_{t\in[0,T]}\left|~\int\limits_0^t\sigma(s,X_s)\,dW_s~\right|~\right].\]
 Then applications of a Krylov estimate, namely Theorem 2.2 in \cite{Zhang_Neu} to the first expectation term and of the inequality of Burkholder, Davis and Gundy (see e.g. \cite{RY05} Corollary IV.4.2) to the second yield
 \[ \mathbb{E}\left[\sup_{t\in[0,T]}|X_t|\right] \leq |x| + C\|b\|_{L_p^q(T)} + C\mathbb{E}\left[\left(\int\limits_0^T\sigma(s,X_s)^2\,ds\right)^\frac{1}{2}\right].\]
 Since $\sigma$ is bounded and $b\in L_p^q(T)$, this is finite. Furthermore,
 \begin{align*}
  \sup_{t\in[0,T]}\mathbb{E}\left[|X_t|^2\right] \leq 2|x|^2 & + 4\sup_{t\in[0,T]}\mathbb{E}\left[~\left|~\int\limits_0^tb(s,X_s)\,ds~\right|^2\right]\\
                                                             & + 4\sup_{t\in[0,T]}\mathbb{E}\left[~\left|~\int\limits_0^t\sigma(s,X_s)\,dW_s~\right|^2\right].
 \end{align*}
 We apply H\"older's inequality to the first expectation and the multidimensional It\^o Isometry to the second one to receive
 \[\sup_{t\in[0,T]}\mathbb{E}\left[|X_t|^2\right]\leq 2|x|^2 + 4T\mathbb{E}\left[~\int\limits_0^T|b(s,X_s)|^2\,ds~\right] + 4\sup_{t\in[0,T]}\mathbb{E}\left[\int\limits_0^t|\sigma(s,X_s)|^2\,ds\right].\]
 Again, we use Theorem 2.2 of \cite{Zhang_Neu} and Assumption \ref{assumptions} (c1), (c4) to obtain that this is finite.
\end{proof}

The following proof of the convergence of the drift term becomes simple with the help of the Krylov estimate Theorem 2.2 of \cite{Zhang_Neu}. The price to pay is the factor two in the assumptions on $p$ and $q$.

\begin{proof}(Proof of Lemma \ref{convergence of the drift})
 Theorem 2.2 of \cite{Zhang_Neu} for $\frac{p}{2}$, $\frac{q}{2}$ and an application of Lemma \ref{Ungleichungen} yields
 \begin{align*}
  \mathbb{E} & \left[\int\limits_0^T\left|b^{(n)}(t,X_t^{(1)})-b^{(n)}(t,X_t^{(2)})\right|^2\,dt\right]\\
             & \leq 2\mathbb{E}\left[\int\limits_0^T\left|b^{(n)}(t,X_t^{(1)})\right|^2\,dt\right] + 2\mathbb{E}\left[\int\limits_0^T\left|b^{(n)}(t,X_t^{(2)})\right|^2\,dt\right]\\
             & \leq C(d,p,q,T,c_\sigma,\tilde{c}_\sigma,\|b\|_{L_p^q(T)})\|b^{(n)}\|_{L_p^q(T)}^2\\
             & \leq C(d,p,q,T,c_\sigma,\tilde{c}_\sigma,\|b\|_{L_p^q(T)})\frac{1}{2^{2(n+1)}}\|b\|_{L_p^q(T)}^2\\
             & \xrightarrow{n\to\infty} 0.
 \end{align*}
\end{proof}

\begin{proof}(Proof of Proposition \ref{uniform bound for A_T})
 Considering $\sigma^{(n)}$ we find that:
 \[\partial_{x_i}\sigma^{(n)} = \left(\partial_{x_i}\partial_xU^{(n)}\right)\sigma + \partial_xU^{(n)}\partial_{x_i}\sigma + \partial_{x_i}\sigma.\]
 We use that $\sigma$ is bounded and $\partial_x\sigma\in L_p^q(T)$, that $\partial_xU^{(n)}$ is uniformly bounded by $\frac{1}{2}$ and $\partial_x^2U^{(n)}$ is equibounded in $L_p^q(T)$ (see Lemma \ref{Ungleichungen}) to deduce that
 \[\|\partial_x\sigma^{(n)}\|_{L_p^q(T)}\leq C \text{ uniformly in } n.\]
 Additionally, $\sigma^{(n)}$ is continuous, since $\partial_xU^{(n)}$ is H\"older continuous. Then there exists a sequence of continuous functions $(u_m)_m$, which are differentiable with respect to $x$ in the ordinary sense, such that $u_m\to\sigma^{(n)}$ uniformly on $[0,T]\times\mathbb{R}^d$
 and
 \[\|\partial_xu_m\|_{L_p^q(T)}\leq\|\partial_x\sigma^{(n)}\|_{L_p^q(T)}~~\forall~ m\in\mathbb{N}.\]
 The existence of such a function can be obtained by mollification. Define $X_t^\lambda:=\lambda X_t^{(1)}+(1-\lambda)X_t^{(2)}$.
 Then we have with Lemma \ref{Ungleichungen} (iv) and uniform convergence
 \begin{align*}
  & \mathbb{E}\left[\exp\left(\int\limits_0^T\frac{\left|\sigma^{(n)}(t,X_t^{(1)})-\sigma^{(n)}(t,X_t^{(2)})\right|^2}{\left|Y_t^{(1,n)}-Y_t^{(2,n)}\right|^2}\textbf{1}_{\{Y_t^{(1,n)}\neq Y_t^{(2,n)}\}}\,dt\right)\right]\\
  & \leq \mathbb{E}\left[\exp\left(4\int\limits_0^T\frac{\left|\sigma^{(n)}(t,X_t^{(1)})-\sigma^{(n)}(t,X_t^{(2)})\right|^2}{\left|X_t^{(1)}-X_t^{(2)}\right|^2}\textbf{1}_{\{X_t^{(1)}\neq X_t^{(2)}\}}\,dt\right)\right].\\
  & = \lim_{m\to\infty}\mathbb{E}\left[\exp\left(4\int\limits_0^T\frac{\left|u_m(t,X_t^{(1)})-u_m(t,X_t^{(2)})\right|^2}{\left|X_t^{(1)}-X_t^{(2)}\right|^2}\textbf{1}_{\{X_t^{(1)}\neq X_t^{(2)}\}}\,dt\right)\right]\intertext{}
  & \leq \lim_{m\to\infty}\mathbb{E}\left[\exp\left(4\int\limits_0^T\int\limits_0^1|\partial_xu_m(t,X_t^\lambda)|^2\,d\lambda \,dt\right)\right]\\
  & \leq \lim_{m\to\infty}\int\limits_0^1\mathbb{E}\left[\exp\left(4\int\limits_0^T|\partial_xu_m(t,X_t^\lambda)|^2\,dt\right)\right]\,d\lambda.
 \end{align*}
 Now, choose $\mu>0$ so small that $(d/p+2/q)(1+\mu)<1/2$ holds.
Then we have for $\beta>0$ with Young's and H\"older's inequality
 \begin{align}
  & \mathbb{E}\left[\exp\left(\int\limits_0^T\frac{\left|\sigma^{(n)}(t,X_t^{(1)})-\sigma^{(n)}(t,X_t^{(2)})\right|^2}{\left|Y_t^{(1,n)}-Y_t^{(2,n)}\right|^2}\textbf{1}_{\{Y_t^{(1,n)}\neq Y_t^{(2,n)}\}}\,dt\right)\right]\notag\\
  & \leq \lim_{m\to\infty}\int\limits_0^1\mathbb{E}\left[\exp\left(\frac{1}{\mu+1}\left(\beta\int\limits_0^T|\partial_xu_m(t,X_t^\lambda)|^2\,dt\right)^{1+\mu}\right.\right.\\
  & ~~~~~~~~~~~~~~~~~~~~~~~~~~~~~~~~~~~~\left.\left.\phantom{\left(\int\limits_0^T\right)^{1+\mu}}+\frac{\mu}{1+\mu}\left(\frac{4}{\beta}\right)^\frac{1+\mu}{\mu}\right)\right]\,d\lambda\notag\\
  & \leq \exp\left(\frac{\mu}{1+\mu}\left(\frac{4}{\beta}\right)^\frac{1+\mu}{\mu}\right)\label{exponential estimate}\\
  & ~~~~~~~~~\cdot\lim_{m\to\infty}\int\limits_0^1\mathbb{E}\left[\exp\left(\int\limits_0^T\frac{\beta^{1+\mu}}{1+\mu}T^\frac{\mu}{1+\mu}|\partial_xu_m(t,X_t^\lambda)|^{2(1+\mu)}\,dt\right)\right]\,d\lambda.\notag
 \end{align}
 Furthermore, we have with Theorem 2.2 from \cite{Zhang_Neu} for all $0\leq t_0\leq T$
 \begin{align*}
  \mathbb{E} & \left[~\int\limits_{t_0}^T\frac{\beta^{1+\mu}}{1+\mu}T^\frac{\mu}{1+\mu}\left|\partial_xu_m(t,X_t^\lambda)\right|^{2(1+\mu)}\,dt~\Bigg|~\mathcal{F}_{t_0}\right]\\
             & \leq C(d,p,q,\mu,T,c_\sigma,\tilde{c}_\sigma,\|b\|_{L_p^q(T)})\frac{\beta^{1+\mu}}{1+\mu}T^\frac{\mu}{1+\mu}\||\partial_xu_m|^{2(1+\mu)}\|_{L_\frac{p}{2(1+\mu)}^\frac{q}{2(1+\mu)}(T)}\\
             & \leq C(d,p,q,\mu,T,c_\sigma,\tilde{c}_\sigma,\|b\|_{L_p^q(T)})\beta^{1+\mu}\|\partial_x\sigma^{(n)}\|_{L_p^q(T)}^{2(1+\mu)}.
 \end{align*}
 Since $\|\partial_x\sigma^{(n)}\|_{L_p^q(T)}$ is equibounded, we can choose $\beta$ so small that this is less than some $0<\alpha<1$ for all $n\in\mathbb{N}$. Then we have by Lemma \ref{Khasminski lemma} and inequality \eqref{exponential estimate} that
 \begin{align*}
  \mathbb{E} \left[e^{A_T^{(n)}}\right] & \leq \exp\left(\frac{\mu}{1+\mu}\left(\frac{4}{\beta}\right)^\frac{1+\mu}{\mu}\right)\frac{1}{1-\alpha}\\
                                                                            & \leq C,
 \end{align*}
 where $C$ does not depend on $n$.
\end{proof}

\section*{Acknowledgement}
The author is very grateful to Michael R\"ockner for useful discussions.


\appendix

\section{}

\begin{proposition}\label{Ito-Formel}(It\^o's formula)
Let $(c1)-(c4)$ of Assumption \ref{assumptions} be fulfilled, $X_t$ a solution to \eqref{SDE} and $u\in W_{q,p}^{1,2}(T)$. Then for $0\leq s\leq t\leq T$ we have
 \begin{align*}
  u(t,X_t) & = u(s,X_s)+\int\limits_s^t\partial_tu(r,X_r)\,dr+\int\limits_s^t\partial_xu(r,X_r)b(r,X_r)\,dr\\
           & ~~~ +\int\limits_s^t\partial_xu(r,X_r)\sigma(r,X_r)\,dW_r\\
           & ~~~ +\frac{1}{2}\int\limits_s^t\sum\limits_{i=1}^d\sum\limits_{j=1}^d(\sigma\sigma^*(r,X_r))_{ij}\partial_{x_ix_j}^2u(r,X_r)\,dr~~~~~~~~\mathbb{P}\text{-almost surely.}
 \end{align*}
\end{proposition}
This result can be obtained by approximation with smooth functions as in \cite[][Theorem 3.7]{KR05}  with the help of \cite[][Theorem 2.2]{Zhang_Neu}.

\begin{lemma}\label{Khasminski lemma}
 Let $f:[0,T]\times\mathbb{R}^d\to\mathbb{R}$ be a nonnegative measurable function and $\gamma$ an arbitrary stopping time. Assume that $X_t$ is an adapted process such there exists a constant $\alpha<1$ with
 \[\textbf{1}_{\{t_0\leq\gamma\}}\mathbb{E}\left[~\int\limits_{t_0}^{T\wedge\gamma}f(t,X_t)\,dt~\Bigg|~\mathcal{F}_{t_0}\right]\leq\alpha~~~~~\mathbb{P}\text{-a.\,s.}~~\forall~0\leq t_0\leq T.\]
 Then we have
 \[\mathbb{E}\left[\exp\left(\int\limits_0^{T\wedge\gamma}f(t,X_t)\,dt\right)\right]\leq\frac{1}{1-\alpha}.\]
\end{lemma}
This is a slightly more general version of Khasminski's Lemma which can be obtained by rewriting the exponential series and using properties of the conditional expectation.

\end{document}